\documentclass[11pt,amssymb,amscd]{amsart}
\usepackage{color}
\usepackage{amssymb}
\usepackage{amsmath}
\usepackage{amsthm} 
\usepackage{mathrsfs}                                                           
\usepackage{slashbox}  
\usepackage[all]{xy}
\usepackage{graphicx}
\usepackage{wrapfig}

\def\2{{1\over 2}}

\def\d{\partial}
\newcommand{\rf}[1]{(\ref{#1})}
\def\b{\bar}

\renewcommand{\t}{\tilde}
\newcommand{\p}{\partial}

\def\p{\partial}

\def\cQ{\mathcal{Q}}
\def\cE{\mathcal{E}}
\def\bcE{\bar{\mathcal{E}}}
\def\cF{\mathcal{F}}

\def\b{\bar}

\def\<{\langle}
\def\>{\rangle}
\def\+{\dagger}

\begin{document}
\title{Beltrami-Courant Differentials and $G_{\infty}$-al\-geb\-ras}
\author{Anton M. Zeitlin}
\address{ 
\newline Department of Mathematics,\newline
Columbia University,
\newline 2990 Broadway, New York,\newline
NY 10027, USA.
\newline
Institut des Hautes \'Etudes Scientifiques,\newline
35 Route de Chartres,\newline 
91440 Bures-sur-Yvette\newline 
IPME RAS, V.O. Bolshoj pr., 61, 199178, \newline
St. Petersburg.\newline
zeitlin@math.columbia.edu\newline
http://math.columbia.edu/$\sim$zeitlin \newline
http://www.ipme.ru/zam.html  }

\begin{abstract}
Using the symmetry properties of two-di\-men\-sional  sigma models, we introduce a notion of the Beltrami-Courant differential, so that there is a natural homotopy Gerstenhaber algebra related to it. 
We conjecture that the generalized Maurer-Cartan equation for the corresponding $L_{\infty}$ subalgebra gives solutions to the Einstein equations. 
\end{abstract}
\maketitle
\section{Introduction}

The geometric and algebraic properties of two-dimensional sigma-models lead to a lot of important discoveries in mathematics. 
One of the most interesting topics, emerged this way in the last decade is the study of gerbes of chiral differential operators, 
which give the proper mathematical description of the simplest first-order sigma-models. In \cite{lmz}, it was shown that the classical actions of 
the  standard second-order sigma-models can be reformulated under certain conditions (one of which is the introduction of complex structure) in terms of 
perturbed first-order ones. In the same article, it was also suggested that the conformal invariance conditions for the perturbed sigma model, 
which have the form of the Einstein (and higher order) equations, will have a homotopical meaning as generalized Maurer-Cartan equations for certain $L_{\infty}$ algebra. In this paper, 
we show that there is a larger structure, namely of homotopy Gerstenhaber algebra, so that the desired $L_{\infty}$ structure is a part of it.

The central object in the construction is the vertex algebroid with a Calabi-Yau structure and its classical limit, 
the Courant algeboid. In \cite{zeit} we associated to every positively graded vertex operator algebra (VOA) the homotopy Gerstenhaber algebra, 
which, according to the work of \cite{lz}, \cite{huangzhao}, \cite{kvz}, \cite{voronov} can be extended to $G_{\infty}$ algebra \cite{tamarkin} and even to $BV_{\infty}$ algebra \cite{gorbounov}, \cite{vallette}. 
The relationship between vertex algebroid and vertex algebra is similar to 
the relationship between Lie algebra and its universal enveloping algebra \cite{gms}.  
We show here that the correspondence constructed in \cite{zeit} can be reformulated by consructing a functor from the category of 
vertex algebroids to the category of $G_{\infty}$-algebras. 
Another important observation of the article \cite{zeit} is that one can construct a quasiclassical limit of the resulting $G_{\infty}$ algebra, so that the 
operations become covariant, i.e. can be expressed via the operations of Courant algebroid only. This $G_{\infty}$ algebra is much easier to 
grasp: its $C_{\infty}$ and $L_{\infty}$ subalgebras reduce to $C_3$ and $L_3$ algebras, where the $L_3$-algebra is the extension of the $L_3$ algebra of Roytenberg and Weinstein \cite{weinstein}. 
An example of the above construction we need in this paper is the $G_{\infty}-$ 
algebra for the vertex algebroid on the space of holomorphic sections of $T^{(1,0)}M\oplus {T^*}^{(1,0)}M$ and its antiholomorphic counterpart, so that the 
corresponding vertex algebra gives (locally) a description of 
the unperturbed first-order sigma-model. The appropriately completed tensor product of corresponding "holomorphic" and "antiholomorphic" homotopy 
Gerstenhaber  algebras gives the homotopy Gerstenhaber algebra and we conjecture that this homotopy Gerstenhaber algebra can be extended to $G_{\infty}$ algebra.  
The Maurer-Cartan elements for the resulting $L_{\infty}$-subalgebra 
are parametrized by the perturbation terms of the first-order sigma model, i.e. by the sections of $\Gamma((T^{(1,0)}M\oplus {T^*}^{(1,0)}M)\otimes (T^{(0,1)}M\oplus {T^*}^{(0,1)}M))\oplus \mathcal{C}(M)$. 
We call the sections from the first summand as $Beltrami-Courant$ $differentials$, justifying that name by its symmetry transformations of the first-order sigma-model, 
which are very similar to the ones of Beltrami differentials on Riemann surfaces and by the fact that the infinitesimal formula is expressed algebraically 
via the operations on Courant algebroid. The sections of the second term in the summand will be called $normalized$ $dilaton$ $fields$. 

It is possible to show that there is a subcomplex in the complex on which the homotopy Gerstenhaber algebra is defined, so that all higher homotopies starting from trilinear ones vanish (in fact, the resulting Gerstenhaber algebra is a BV-algebra \cite{getzler}). 
We show that the Maurer-Cartan equation of the corresponding differential graded Lie algebra is equivalent to Einstein Equations with dilaton and B-field, 
if the bivector field from $\Gamma(T^{(1,0)}M\otimes T^{(0,1)}M)$, which parametrizes the Maurer-Cartan element, gives rise to the Hermitian metric.
This leads to the conjecture, in view of the relation between first and second-order sigma models that the generalized Maurer-Cartan equations 
(GMC) for $L_{\infty}$-algebra on the full complex give Einstein equations with B-field and dilaton, parametrized by the Beltrami-Courant differential. 
We justify the conjecture by showing that the symmetries of GMC reproduce the infinitesimal diffeomorphisms and gauge transformations of a B-field up 
to the second order in the Beltrami-Courant differential.

The structure of the paper is as follows. 
In Section 2 we study the classic action functionals for first- and second-order sigma models and the relationship between their symmetries. 
This leads to the definition of the Beltrami-Courant differential and its symmetries, e.g. under diffeomorphism transformations. 
In Section 3, we discuss Vertex/Courant algebroids with the Calabi-Yau structure,  
related $G_{\infty}$-algebras and their classical limits. 
In Section 4 we describe the relation of these algebras to 
Einstein equations with B-field and dilaton, parametrized by Beltrami differential. First we describe the simplest case, 
when the $G_{\infty}$ algebra is reduced to the Gerstenhaber algebra, then we formulate the conjecture regarding more general Einstein 
equations and support it by calculation of the symmetry transformations.\\
  
\noindent{\bf Acknowledgements.} The author is indebted to M. Kontsevich, S. Merkulov and B. Valette for fruitful discussions and remarks, and to A.N. Fedorova for the careful reading of the manusript.

\section{Sigma-models and Beltrami-Courant Differentials}
In this section, we introduce the first object of interest: Beltrami-Courant differential. 
We derive its definition from the symmetries of the classical sigma-model actions. Let us consider a complex Riemann surface $\Sigma$, a complex manifold $M$ of dimension $d$, and a map $X:\Sigma\to M$. Then one can write the following action functional:
\begin{eqnarray}\label{free}
S_0=\frac{1}{2\pi ih}\int_\Sigma (\langle p\wedge\bar{\partial} X\rangle-
\langle \bar{p}\wedge{\partial} X\rangle),
\end{eqnarray}
where $p$ and $\bar{p}$ belong to $X^*(\Omega^{(1,0)}(M))\otimes \Omega^{(1,0)}(\Sigma)$ and 
$X^*(\Omega^{(0,1)}(M))\otimes \Omega^{(0,1)}(\Sigma)$ correspondingly and 
$\langle\cdot, \cdot\rangle$ stands for standard pairing. 
In the following $X^i, X^{\bar j}$ stand for the pull-backs of the coordinate functions on $M$ with respect to $X$. 
This action has the following symmetries (we write them in components in the infinitesimal form):
\begin{eqnarray}\label{set1h}
&&X^i\to X^i-v^i(X), \quad p_i\to p_i+\p_i v^k p_k,\nonumber\\
&& X^{\b i}\to X^{\b i}-v^{\b i}(\b X), \quad p_{\b i}\to p_{\b i}+\p_{\b i} v^{\b k} p_{\b k} .
\end{eqnarray}
Here, the generators of the infinitesimal transformations $v, \b v$ are the elements of $\Gamma(\mathcal{O}(T^{(1,0)}M))$ and 
$\Gamma(\mathcal{O}(T^{(0,1)}M))$ correspondingly, i.e. $v^i$ ($v^{\b i}$) 
are (anti)holomorphic. These symmetries illustrate invariance under the holomorphic coordinate transformations.  
There is another set of symmetries, induced by the (anti)holomorphic 1-forms. Let 
$\omega\in \Gamma(\mathcal{O}({T^*}^{(1,0)}M))$ and $\b\omega\in \Gamma(\bar{\mathcal{O}}({T^*}^{(0,1)}M))$. Then the action \rf{free} is invariant 
under the transformatiosn of $p, \b p$:
\begin{eqnarray}\label{set2h}
p_i\to p_i-\p X^k(\p_k\omega_i-\p_i\omega_k), \quad p_{\b i}\to p_{\b i}-\bar{\p} X^{\b k}(\p_{\b k}\omega_{\b i}-\p_{\b i}\omega_{\b k}).
\end{eqnarray}
We want to generalize the action \rf{free} so that it would be invariant under the diffeomorphism transformations and nonholomorphic generalizations of \rf{set2h}. 
In order to do that, one has to introduce extra (perturbation) terms to the action \rf{free}. Let us see how it works with an example. Suppose $v^i, v^{\b i}$ in the 
formulas \rf{set1h} are not holomorphic anymore, then $S_0$ won't be invariant and there will be an extra contribution to $S_0$: 
\begin{eqnarray}
 \delta S_0=-\frac{1}{2\pi ih}\int_{\Sigma} (\langle \b \p v, p\wedge \b \p X\rangle +\langle \p \b v , \b p\wedge \p X\rangle).
\end{eqnarray}
Therefore, to compensate this term, it makes sense to add extra terms to the action of the form
\begin{eqnarray}
 \delta S_{\mu}=-\frac{1}{2\pi ih}\int_{\Sigma} (\langle \mu, p\wedge \b \p X\rangle +\langle \b \mu, \p X\wedge \b p\rangle),
\end{eqnarray}
where $\mu\in \Gamma (T^{(1,0)}M\otimes {T^*}^{(0,1)}(M))$, $\bar{\mu}\in \Gamma (T^{(0,1)}M\otimes {T^*}^{(1,0)}(M))$, so that upon the $(v, \b v)$ transformations 
$\mu, \bar{\mu}$ should be modified as follows:
\begin{eqnarray}
\mu\to \mu-\b \p v+\dots, \quad \b \mu\to \b \mu-\p \b v+\dots,
\end{eqnarray}
where dots stand for terms higher in $\mu$ and $\bar \mu $. Continuing and further applying this approach to the non(anti)holomorphic generalizations of the transformations 
\rf{set1h}, \rf{set2h} we find that we have to add the following terms to $S_0$, such that the resulting action is: 
\begin{eqnarray}
&&\tilde{S}=
\frac{1}{2\pi ih}\int_\Sigma (\langle p\wedge\bar{\partial} X\rangle-
\langle \bar{p}\wedge{\partial} X\rangle -\\
&&\langle \mu, p\wedge \b \p X\rangle -\langle \b \mu, \p X\wedge \b p\rangle-\langle b, \p X\wedge \b \p X\rangle),\nonumber
\end{eqnarray}
where $b\in \Gamma ({T^*}^{(1,0)}M\otimes {T^*}^{(0,1)}M)$.
The resulting symmetry transformations generated by $(v, \b v)$ can be written as follows:
\begin{eqnarray}\label{mutransf}
 &&\mu^{i}_{\bar{j}} \rightarrow \mu^{i}_{\bar{j}} -
\p_{\bar{j}}v^i + v^{k}\p_k\mu^{i}_{\bar{j}} +
v^{\bar{k}}\p_{\bar{k}}\mu^{i}_{\bar{j}}+
\mu^{i}_{\bar{k}}\p_{\bar{j}}v^{\bar{k}} -
\mu^k_{\bar{j}}\p_kv^i
 + \mu^i_{\bar{l}}\mu^k_{\bar{j}}\p_k v^{\bar{l}},\\ 
&&b_{i{\bar j}} \rightarrow b_{i{\bar j}} + v^k\p_k b_{i{\bar j}} + v^{\bar{k}}\p_{\bar{k}} b_{i{\bar j}}
+ b_{i{\bar k}}\p_{\bar{j}}v^{\bar{k}}+b_{l{\bar j}}\p_i v^l+b_{i{\bar k}}\mu^{k}_{\bar j}\p_kv^{\bar{k}}
+b_{l{\bar j}}{\bar\mu}^{\bar k}_i\p_{\bar k}v^l,\nonumber
 \end{eqnarray}
 and the formula for the transformation of $\bar{\mu}$ can be obtained from the one of $\mu$ by formal complex conjugation. 
This leads to the symmetry of the action $\tilde S$ if 
\begin{eqnarray}
&&X^i\to X^i-v^i(X, \b X), \quad p_{i} \rightarrow p_{i} + p_k \p_iv^k - p_k\mu^k_{\bar{l}}\p_iv^{\bar{l}}
- b_{j{\bar k}}\p_i v^{\bar k}\p X^j,\\
&& X^{\b i}\to X^{\b i}-v^{\b i}(X,\b X), \quad \bar{p}_{\bar{i}} \rightarrow \bar{p}_{\bar{i}} + \bar{p}_{\bar k} \p_{\bar i}v^{\bar {k}} - {\bar p}_{\bar{k}}
{\bar \mu}^{\bar k}_{l}\p_iv^{l}
- b_{\bar{j}k}\p_{\bar{i}} v^{k}\bar{\p} X^{\bar {j}}.\nonumber
\end{eqnarray}
Therefore, the resulting action is invariant under the action of the infinitesimal diffeomorphism group. The component formulas \rf{mutransf} were first discovered in \cite{gamayun}.  
Similarly, we obtain that the transformations 
\begin{eqnarray}\label{omtransf}
&&b_{i\bar{j}}\to  b_{i\bar{j}}+\p_{\bar j}\omega_i-\p_i\omega_{\bar j}+\mu^i_{\bar j}(\p_i\omega_k-\p_k\omega_i)+\nonumber\\
&&\bar{\mu}^{\bar s} _i(\p_{\bar j}\omega_{\bar s}-
\p_{\bar s}\omega_{\bar j})+
{\bar \mu}^{\bar i}_j\mu_{\bar k}^s(\p_s\omega_{\bar i}-\p_{\bar i}\omega_s)
\end{eqnarray}
accompanied with 
\begin{eqnarray}
&& p_i\to p_i-\p X^k(\p_k\omega_i-\p_i\omega_k)-\p_{\b r}\omega_i\p X^{\bar r}-{\bar \mu}^{\bar s}_k\p_i\omega_{\bar s}\p X^k, \nonumber\\
&& p_{\b i}\to p_{\b i}-{\bar\p} X^{\b k}(\p_{\b k}\omega_{\b i}-{\p}_{\b i}\omega_{\b k})-\p_{r}\omega_{\bar i}{\bar\p} X^{r}-{\mu}^{s}_{\bar k}\p_i\omega_{s}{\bar\p} X^{\bar k}.
\end{eqnarray}
leave $\tilde{S}$ invariant. 
Hereinafter, it is useful to consider $\mu, {\bar \mu}, b$ as matrix elements of $\tilde{\mathbb{M}}\in \Gamma((T^{(1,0)}M\oplus {T^*}^{(1,0)}M)\otimes(T^{(0,1)}M\oplus {T^*}^{(0,1)}M))$, i.e. 
\begin{equation}
\tilde{\mathbb{M}}=\begin{pmatrix} 0 & \mu \\ 
\bar{\mu} & b \end{pmatrix}.
\end{equation}
For simplicity of notation let us define $E=TM\oplus T^*M$, also 
$\mathcal{E}=T^{(1,0)}M\oplus {T^*}^{(1,0)}M$ and $\bar{\mathcal{E}}=T^{(0,1)}M\oplus {T^*}^{(0,1)}M$, so that $E=\mathcal{E}\oplus\bar{\mathcal{E}}$.

Let $\alpha\in \Gamma(E)$, i.e. $\alpha=(v, \bar v, \omega, \bar \omega)$, where  $v, \b v$ are the elements of $\Gamma(T^{(1,0)}M)$ and 
$\Gamma(T^{(0,1)}M)$ correspondingly 
and $\omega\in \Omega^{(1,0)}(M)$, $\b\omega\in \Omega^{(0,1)}(M)$. Next, we introduce operator the operator 
$D:\Gamma(E)\to \Gamma(\cE\otimes\bcE)$, such that

\[ D\alpha=\left( \begin{array}{cc}
0 & {\bar \p }v\\
{\p \bar v} & \p{\bar\omega}-{\bar\p} \omega \end{array} \right).\]

Then the transformation of $\tilde{\mathbb{M}}$ under \rf{mutransf}, \rf{omtransf} can be expressed by the following formula:
\begin{eqnarray}\label{algsym}
\tilde{\mathbb{M}}\to \tilde{\mathbb{M}}-D\alpha+ \phi_1(\alpha,\tilde{\mathbb{M}})+\phi_2(\alpha, \tilde{\mathbb{M}},\tilde{\mathbb{M}}),
\end{eqnarray}
where we separated terms with linear and 
bilinear dependence on $\tilde M$, denoting them by $\phi_1$ and $\phi_2$ correspondingly.   
There is a hidden algebraic meaning of $\phi_1$ and $\phi_2$ operations on $\alpha$ and $\tilde{M}$. In order to uncover it, we have to use jet bundles.   
Namely, let us consider 
\begin{eqnarray}\label{jets}
&&\xi\in J^{\infty}(\mathcal{O}_M)\otimes 
J^{\infty}({\bar{\mathcal{O}}}(\bcE))
\oplus J^{\infty}({\mathcal{O}}(\cE))
\otimes J^{\infty}({\bar{\mathcal{O}}}_M),\\ 
&&\mathbb{L}\in J^{\infty}({\mathcal{O}}(\cE))\otimes J^{\infty}({\bar{\mathcal{O}}}(\bcE)),\nonumber
\end{eqnarray}
where $J^{\infty}(E)$, for any bundle $E$ over $M$ stands for the corresponding $\infty$-jet bundle of $E$.  
In other words, let 
\begin{eqnarray}\label{algnot}
&&\xi=\sum_Jf^J\otimes {\bar{b}}^J+\sum_Kb^K\otimes {\bar{f}}^K,\nonumber\\ 
&&\mathbb{L}=\sum_I a^I\otimes \bar{a}^I, 
\end{eqnarray}
where $a^I, b^J\in J^{\infty}({\mathcal{O}}(\cE))$, $f^I\in 
J^{\infty}(\mathcal{O}_M)$ and ${\bar{a}}^I, {\bar{b}}^J\in J^{\infty}({\bar{\mathcal{O}}}(\bcE))$, 
${\bar{f}}^I\in 
J^{\infty}(\bar{\mathcal{O}}_M)$. Then we can introduce the operation $\phi_1(\xi, \mathbb{L})$ as follows:
\begin{eqnarray}\label{bialgop}
\phi_1(\xi,\mathbb{L})=
\sum_{I,J}[b^J, a^I]_D\otimes 
{\bar{f}}^{J}{\bar{a}}^I+\sum_{I,K}f^Ka^I\otimes[{\bar{b}}^K, {\bar{a}}^I]_D,
\end{eqnarray}
where $[\cdot, \cdot]_D$ is a Dorfman bracket, see e.g. \cite{roytenberg1} or the next section. Completing the tensor products 
in \rf{jets}, one can introduce the operation on $\alpha\in \Gamma(E)$ and 
$\tilde{\mathbb{M}}\in\Gamma(\cE\otimes \bcE)$, which we also denote as $\phi_1$. 
One can explicitly check that 
\rf{bialgop} leads to the part of \rf{mutransf} and \rf{omtransf}, linear in $\alpha$ and $\tilde{\mathbb{M}}$. The last part, bilinear in 
$\tilde{\mathbb{M}}$, also has an algebraic meaning of a similar kind: returning back to the notation \rf{algnot}, we find that on the jet counterparts 
of $\alpha, \tilde{\mathbb{M}}$, i.e. on $\xi, \mathbb{L}$ the expression for $\phi_2$ is:
\begin{eqnarray}\label{trialgop}
&&\phi_2(\xi, \mathbb{L},\mathbb{L})=\nonumber\\
&&\frac{1}{2}\sum_{I,J,K}\langle b^I, a^K\rangle a^J\otimes \bar{a}^J(\bar{f}^I) {\bar{a}}^K+\frac{1}{2}\sum_{I,J,K}
 a^J(f^I) {a}^K\otimes \langle {\bar{b}}^I, {\bar{a}}^K\rangle {\bar{a}}^J,
\end{eqnarray}
where $\bar{a}^J(\bar{f}^I)$, $a^J(f^I)$ correspond to the action of the differential operator, associated to the vector field, on a function 
($\bar{a}^J(\bar{f}^I)$, $a^J(f^I)$ are set to be zero if $\bar{a}^J$, 
$a^J$ are 1-forms).  
At the same time, the operation $\phi_2$ 
has the following simple description:
\begin{eqnarray}\label{mudamu}
\phi_2(\alpha, \tilde{\mathbb{M}},\tilde{\mathbb{M}})=\tilde{\mathbb{M}}\cdot D\alpha\cdot \tilde{\mathbb{M}}
\end{eqnarray}
if we consider $\tilde{\mathbb{M}}$ as an element of $End(\Gamma(E))$.

Let us notice that we could generalize $\tilde{\mathbb{M}}$ in the following way: in the matrix expression for 
$\tilde{\mathbb{M}}$ let us fill in the empty spot, i.e. let us add extra element $g\in\Gamma(T^{(1,0)}M\otimes T^{(0,1)}M)$. Then the modified 
$\tilde{\mathbb{M}}$, i.e.    
$\mathbb{M}\in\Gamma(\cE\otimes\bcE)$ can be expressed as follows: 
\begin{equation}\label{mmat}
\mathbb{M}=\begin{pmatrix} g & \mu \\ 
\bar{\mu} & b \end{pmatrix}.
\end{equation}
The corresponding action functional is:
\begin{eqnarray}
&&S_{fo}=\frac{1}{2\pi i h}\int_\Sigma (\langle p\wedge\bar{\partial} X\rangle+
\langle \bar{p}\wedge{\partial} X\rangle-\nonumber\\
&& -\langle g, p\wedge \bar{p} \rangle-\langle \mu, p\wedge \b \p X\rangle -\langle \b \mu, \b p\wedge \p X\rangle-\langle b, \p X\wedge \b \p X\rangle).
\end{eqnarray}
It turns out that the symmetries of this action functional can be described by the same formula \rf{algsym}, where algebraic meaning of the 
operations on the jet level is given by the same formulas \rf{bialgop}, \rf{trialgop}, and the formula \rf{mudamu} is also valid. In Appendix, one can 
find the explicit component formulas for the infinitesimal symmetries of the action $S_{fo}$.
The reason for introducing the $g$-term in the action functional is as follows. If the matrix $\{g^{i\bar{j}}\}$ is invertible, 
then using elementary variational calculus, one can find that 
the critical points for $S_{fo}$ are the same as for the 
second-order action functional:
\begin{eqnarray}
&&S_{so}=\\
&&\frac{1}{2\pi h}\int_{\Sigma} d^2 z (g_{i\bar{j}}(\bar{\partial} X^i-\mu^{i}_{\bar{k}}
\bar{\partial} X^{\bar{k}})(\partial X^{\bar{j}}-\bar{\mu}_{k}^{\bar{j}}
\partial X^{k})-b_{i\bar{j}}\partial X^{i}\bar{\partial} X^{\bar{j}}),\nonumber
\end{eqnarray}
which can be re-expressed as
\begin{eqnarray}
S_{so}=\frac{1}{4\pi h}\int_{\Sigma} d^2 z
(G_{\mu\nu}+B_{\mu\nu})\partial X^{\mu}\bar{\partial}X^{\nu},
\end{eqnarray}
where $G$ is a symmetric tensor and $B$ is antisymmetric, indices $\mu,\nu$ run through the set $\{i, \bar{j}\}$. The expression for 
$G$ and $B$ via $\mathbb{M}$ is given by:
\begin{eqnarray}\label{GB}
\label{phytwi}
G_{s\bar{k}}&=&g_{\bar{i}j}
\bar{\mu}^{\bar{i}}_s\mu^{j}_{\bar{k}}+g_{s\bar{k}}-
b_{s\bar{k}}, \quad
B_{s\bar{k}}=g_{\bar{i}j}\bar{\mu}^{\bar{i}}_s\mu^{j}_{\bar{k}}-g_{s\bar{k}}-
b_{s\bar{k}},\\
G_{si}&=&-g_{i\bar{j}}\bar{\mu}^{\bar{j}}_s-g_{s\bar{j}}\bar{\mu}^{\bar{j}}_i
, \quad
G_{\bar{s}\bar{i}}=-g_{\bar{s}j}\mu^{j}_{\bar{i}}-g_{\bar{i}j}\mu^{j}_{\bar{s}},
\nonumber\\
B_{si}&=&g_{s\bar{j}}\bar{\mu}^{\bar{j}}_i-g_{i\bar{j}}\bar{\mu}^{\bar{j}}_s,
\quad
B_{\bar{s}\bar{i}}=g_{\bar{i}j}\mu^{j}_{\bar{s}}-g_{\bar{s}j}\mu^{j}_{\bar{i}},
\nonumber
\end{eqnarray}
where $\{g_{i\bar{j}}\}$ stands for the inverse matrix of $\{g^{i\bar{j}}\}$. Such parametrization of the second-order action in the case when $M$ is a Riemann surface was first introduced in \cite{pz}, \cite{zeitlin}.
The symmetries of the action functional $S_{fo}$ transform 
into infinitesimal diffeomorphism transformations and the 2-form $B$ symmetry
\begin{eqnarray}
&& G\to G-L_{\bf v}G,\quad B\to B-L_{\bf v}B,\\ 
&& B\to B-2d{\bf\boldsymbol \omega},\nonumber
\end{eqnarray}
 if $\alpha=({\bf v}, {\boldsymbol \omega})$, so that ${\bf v}\in \Gamma(TM)$, ${\boldsymbol \omega}\in \Omega^{1}(M)$, i.e. the symmetries of $S_{so}$. 

Let us formulate this as a theorem.\\

\noindent {\bf Theorem 1.1.} {\it Let $\mathbb{M}\in\Gamma(\cE\otimes\bcE)$, parametrized as 
in \rf{mmat}, so that its $\Gamma(T^{(1,0)}M\otimes T^{(0,1)}M)$ part is given by $\{g^{i\bar{j}}\}$, which is invertible, 
then the infinitesimal diffeomorphism transformations of the resulting symmetric and antisymmetric tensors $G$ and $B$ (see \rf{GB}), as well as the 
B-tensor shift by exact 2-forms are encoded in the formula 
\begin{eqnarray}\label{tr}
\mathbb{M}\to \mathbb{M}-D\alpha+ \phi_1(\alpha,\mathbb{M})+\phi_2(\alpha, \mathbb{M},\mathbb{M}),
\end{eqnarray}
where $\alpha\in \Gamma(E)$ and operations $\phi_1,\phi_2$ are defined above.}\\

Note that if $\{G_{\mu\nu}\}$ is invertible and real, it gives rise to the metric tensor. Therefore, since $\mathbb{M}$ parametrizes both 
$G$ and $B$, and transforms according to \rf{tr} under diffeomorphisms, it is analogous to Beltrami differential on the  Riemann surface. 
So, from now on we will call the elements of $\Gamma(\cE\otimes \bcE)$ as Beltrami-Courant differentials, since, as we see in the following sections, they are described by means of the Courant algebroid \cite{courant} structure on $\mathcal{E}, \bcE$. 

\section{Vertex algebroids, $G_{\infty}$-algebra and quasiclassical limit}
In this section, we describe the constructions of the article \cite{zeit} with some modifications and refer the reader to this article for some of the details.

Each of the terms in the classical action $S_0$ 
from which we started the previous section, leads to the quantum theory which is well described locally on open neighborhoods of $M$ by means of vertex algebra generated by operator products
\begin{eqnarray}
X^i(z)p_j(w)\sim \frac{h\delta^i_j}{z-w}, \quad 
X^{\b i}(\b z)p_{\b j}(\b w)\sim \frac{h\delta^{\b i}_{\b j}}{\b z-\b w}
\end{eqnarray}
and globally by means of gerbes of chiral differential operators on $M$ \cite{nekrasov}, \cite{gms}. 
Each of the corresponding vertex algebras, which provide the local description, 
form a $\mathbb{Z}_{+}$-graded vector space $V=\sum_{n=0}^{+\infty}V_n$, so that $V_n$ is determined (see \cite{gms}) 
by means of a vertex algebroid. 
In our case, the vertex algebroid is described by means of the sheaf 
$\mathcal{V}=\mathcal{O}(\cE)\otimes \mathbb{C}[h]\equiv \mathcal{O}(\cE)^h$ (resp. $\bar{\mathcal{O}}(\bcE)^h$), of vector spaces $V_1$, as well as the sheaf of vector spaces $V_0$, which coincides with the structure sheaf $\mathcal{O}_M\otimes \mathbb{C}[h]=\mathcal{O}_M^h$ (resp. $\bar{\mathcal{O}}_M^h$), with certain algebraic operations between them.
 
Let us define a vertex alebroid (see e.g. \cite{gms}, \cite{bressler}) and then study our concrete case in detail. 

A {\em vertex $\mathcal{O}^h_M$-algebroid} is a sheaf of $\mathbb{C}[h]$-vector 
spaces $\mathcal{V}$ with a pairing $
\mathcal{O}_M\otimes_{\mathbb{C}[h]}\mathcal{V}  \to  \mathcal{V}$, i.e.  
$f\otimes v  \mapsto  f*v$ 
such that $1* v = v$, equipped with
a structure of a Leibniz $\mathbb{C}[h]$-algebra 
$[\ ,\ ] :
\mathcal{V}\otimes_{\mathbb{C}[h]}\mathcal{V}\to \mathcal{V}$, a $\mathbb{C}[h]$-linear map of Leibniz algebras $\pi : \mathcal{V}\to \Gamma({TM})$, which 
usually is referred to as an anchor, 
a symmetric $\mathbb{C}[h]$-bilinear pairing $\langle\ ,\ \rangle :
\mathcal{V}\otimes_{\mathbb{C}[h]}\mathcal{V}\to \mathcal{O}_M^h$
a $\mathbb{C}$-linear map $\p : \mathcal{O}_M\to \mathcal{V}$ 
such that
$\pi\circ\partial = 0$,
which satisfy the relations
\begin{eqnarray}
&& f*(g*v) - (fg)*v  =  \pi(v)(f)*\partial(g) +
\pi(v)(g)*\partial(f),\nonumber\\
&&[v_1,f*v_2]  =  \pi(v_1)(f)*v_2 + f*[v_1,v_2], 
\nonumber\\
&&[v_1,v_2] + [v_2,v_1]  =  \partial(\langle v_1,v_2\rangle),
\quad
\pi(f*v) = f\pi(v),  \nonumber\\
&&\langle f*v_1, v_2\rangle  =  f\langle v_1,v_2\rangle -
\pi(v_1)(\pi(v_2)(f)), \nonumber\\
&&\pi(v)(\langle v_1, v_2\rangle)  =  \langle[v,v_1],v_2\rangle +
\langle v_1,[v,v_2]\rangle,\nonumber \\
&&\partial(fg)  =  f*\partial(g) + g*\partial(f), \nonumber\\
&&[v,\partial(f)] =  \partial(\pi(v)(f)), \quad
\langle v,\partial(f)\rangle  =  \pi(v)(f),
\end{eqnarray}
where $v,v_1,v_2\in\mathcal{V}$, $f,g\in\mathcal{O}_M^h$. \\

The correspondence between vertex algebroid and the vertex algebra on each neighborhood $U$ is similar to the correspondence between Lie algebra and its universal enveloping algebra: for more details see \cite{gms}. 

Let us concentrate on the case when 
$\mathcal{V}=\mathcal{O}(\cE)^h$. Explicitly, if 
$f\in \mathcal{O}_M$, $v, v_1, v_2\in \mathcal{O}(T^{(1,0)}M)$, 
$\omega, \omega_1, \omega_2\in \mathcal{O}({T^*}^{(1,0)}M)$, then 
locally in the neighborhood with the coordinates $\{X^i\}$
\begin{eqnarray}
&&\p f=df,\quad \pi(v)f=-hv(f), \quad\pi(\omega)=0,\nonumber\\
&&f*v=fv+hdX^i\p_i\p_jfv^j, \quad f*\omega=f\omega,\nonumber\\
&&[v_1,v_2]=-h[v_1,v_2]_D-h^2dX^i\p_i\p_kv^s_1\p_sv^k_2,\nonumber\\
&&[v,\omega]=-h[v,\omega]_D, \quad 
[\omega, v]=-h[\omega,v]_D, \quad [\omega_1,\omega_2]=0,\nonumber\\  
&&\langle v, \omega\rangle=-h\langle v, \omega\rangle^s, \quad \langle v_1, v_2\rangle=-h^2\p_iv_1^j\p_jv_2^i,\quad \langle\omega_1, \omega_2\rangle=0,
\end{eqnarray}
where $\langle\cdot,\cdot \rangle^s$ is a standard pairing on $\cE$ and $[\cdot,\cdot]_D$ is the Dorfman bracket:
\begin{eqnarray}
&&[v_1,v_2]_D=[v_1,v_2]^{Lie}, \quad [v,\omega]_D=L_v\omega,\nonumber\\
&&[\omega , v]_D=-i_vd\omega,\quad [\omega_1,\omega_2]_D=0.
\end{eqnarray}
In \cite{zeit}, it was shown that given a holomorphic volume form on the open neighborhood $U$ of $M$, one can associate a homotopy Gerstenhaber algebra to the vertex algebroid on $U$  
(although the main emphasis of \cite{zeit} was on $C_{\infty}$ part of 
it). 
This was done by considering semi-infinite complex associated to the vertex algebra: due to the results of\cite{lz}, \cite{kvz}, \cite{huangzhao}, \cite{voronov}, there is a structure of $G_{\infty}$ algebra attached to it if the  central charge of the corresponding Virasoro algebra is 26. Using this fact and considering the subcomplex corresponding to the elements of total conformal weight zero, we find out that the central charge condition can be dropped. 
The resulting complex $(\mathcal{F}^{\cdot}, Q)$ appears to be much shorter that original semi-infinite one, namely it is of the form 
\begin{equation}
0\to \mathcal{F}^0\xrightarrow{Q} \mathcal{F}^1\xrightarrow{Q} \mathcal{F}^2\xrightarrow{Q} \mathcal{F}^3\to 0
\end{equation}
and the action of the differential is defined by means of the following diagram
\begin{eqnarray}
\xymatrixcolsep{30pt}
\xymatrixrowsep{3pt}
\xymatrix{
& \mathcal{V}\ar[ddddr] & \mathcal{V}\ar[ddddr]^{\frac{1}{2}h\cdot{\rm div}}& \\
&& \p &&\\
& \bigoplus & \bigoplus & \\
&& {-\frac{1}{2}h\cdot {\rm div}} &&\\
\mathcal{O}_M^{h}\ar[uuuur]^{\p} & \mathcal{O}_M^{h}\ar[uuuur]\ar[r]_{i\rm{d}} & \mathcal{O}_M^{h}  & \mathcal{O}_M^{h}.
}
\end{eqnarray}
Here $\mathcal{F}^0\cong\mathcal{O}_M^{h}\cong\mathcal{F}_h^3$, $\mathcal{F}_h^1\cong\mathcal{O}_M^{h}\oplus \mathcal{V}\cong\mathcal{F}_h^2$, missing arrows correspond to the zero action of $\mathcal{Q}$ and ${\rm div}$ stands for  
divergence operator with respect to the nonvanishing volume form applied to sections of $\Gamma(U,T^{(1,0)}(M))$. Appropriate analogue of operator ${\rm div}$ in the case of general vertex algebroid is called $Calabi-Yau$ $structure$ on vertex algebroid \cite{gms} (since e.g. in our case to be defined globally $M$ should possess a nonvanishing holomorphic volume form).  
According to \cite{zeit}, this complex has a bilinear operation, which satisfies the Leibniz identity with respect to $Q$, it 
is also homotopy commutative and associative, and can be described by the following table:
\begin{equation}\label{asstab}
(a_1,a_2)_h=
\end{equation}
\begin{equation}
\begin{tabular}{|l|c|c|c|c|c|r|}
\hline
 \backslashbox{$a_2$}{$a_1$}&          $u_1$ & $A_1$ & $v_1$ &$\t A_1$ &$\t v_1$& $\t u_1$ \\
\hline
$u_2$ &                  $u_1u_2$    &$A_1u_2$ &$v_1u_2$&
$\t A_1u_2$&$\t v_1 u_2$ &$\t u_1 u_2$\\
&  & $+\pi(A_1)u_2$  &  &   && \\
\hline
$A_2$     &  $u_1A_2$ & $-[A_1,A_2]+$ &$-v_1A_2+$  
&  $\frac{1}{2}\langle\t A_1,A_2\rangle$&$-\pi(A_2)\t v_1$&0\\                    
& & $\frac{1}{2}\langle A_1, A_2\rangle$  & $\pi(\t v_1)A_2$  &   && \\
\hline
$v_2$& $u_1\t u_2$ &  $A_1v_2$ & 0 & $-\pi(\t A_1) v_2 $ &$-\t v_1v_2$& 0\\
\hline
$\t A_2$ & $u_1\t A_2$ &$\frac{1}{2}\langle A_1,\t A_2\rangle$&  
$-\pi(\t A_2)v_1$ & 0& 0&0\\
\hline
$\t v_2$& $u_1\t u_2$ &   $\pi(A_1)\t v_2$& $-v_1\t v_2$ & 0  &0& 0\\
\hline
$\t u_2$    & $u_1\t u_2$ &   0& 0 & 0  &0& 0\\
\hline
\end{tabular},\nonumber
\end{equation}
\vspace{3mm}
where $u_i\in\mathcal{F}_h^0$, $(v_i,A_i)\in\mathcal{F}_h^1$, $(\t v_i,\t A_i)\in\mathcal{F}_h^2$, $\t u_i\in \mathcal{F}_h^3$. 

We note that there is an operator $\mathbf{b}$ of degree -1 on 
$(\mathcal{F}_h^{\cdot}, Q)$ which anticommutes with $Q$:  
\begin{eqnarray}
\xymatrixcolsep{30pt}
\xymatrixrowsep{3pt}
\xymatrix{
& \mathcal{V} & \mathcal{V}\ar[l]_{-i\rm{d}}& \\
& \bigoplus & \bigoplus & \\
\mathcal{O}_M^{h} & \mathcal{O}_M^{h}\ar[l]_{i\rm{d}} & \mathcal{O}_M^{h}& \mathcal{O}_M^{h}\ar[l]_{-i\rm{d}}
} 
\end{eqnarray}
This operator gives rise to the bracket operation 
\begin{eqnarray}\label{brack}
(-1)^{|a_1|}\{a_1,a_2\}_h=\mathbf{b}(a_1,a_2)_h-(\mathbf{b}a_1,a_2)_h-(-1)^{|a_1|}(a_1\mathbf{b}a_2)_h,
\end{eqnarray}
which satisfies quadratic relations together with $(\cdot, \cdot)_h$ and $Q$, which follows from the properties of the vertex algebra \cite{lz}. On the cohomology of $Q$ these relations turn 
into defining properties of Gerstenhaber algebra. Namely, the following Proposition holds.\\

\noindent {\bf Proposition 3.1.}\cite{zeit} {\it Symmetrized versions of operations \rf{asstab} together with \rf{brack} satisfy the relations of the homotopy Gerstenhaber algebra, which follows from these relations: 
\begin{eqnarray}\label{lzrel}
&&Q(a_1,a_2)_h=(Q a_1,a_2)_h+(-1)^{|a_1|}(a_1,Q a_2)_h,\\
&&(a_1,a_2)_h-(-1)^{|a_1||a_2|}(a_2,a_1)_h=\nonumber\\
&&Qm(a_1,a_2)+m(Qa_1,a_2)+(-1)^{|a_1|}m(a_1,Qa_2),\nonumber\\
&& Qn_h(a_1,a_2,a_3)_h+n_h(Qa_1,a_2,a_3)+(-1)^{|a_1|}n_h(a_1,Qa_2,a_3)+\nonumber\\
&&(-1)^{|a_1|+|a_2|}n_h(a_1,a_2,Qa_3)=((a_1,a_2)_h,a_3)_h-(a_1,(a_2,a_3)_h)_h,\nonumber\\
&&\{a_1,a_2\}+(-1)^{(|a_1|-1)(|a_2|-1)}\{a_2,a_1\}=\nonumber\\
&&(-1)^{|a_1|-1}(Qm_h'(a_1,a_2)-m_h'(Qa_1,a_2)-(-1)^{|a_2|}m_h'(a_1,Qa_2)),
\nonumber\\
&& \{a_1,(a_2,a_3)_h\}_h=(\{a_1,a_2\}_h,a_3)_h+(-1)^{(|a_1|-1)||a_2|}(a_2,\{a_1, a_3\}_h)_h,\nonumber\\
&&\{(a_1,a_2)_h,a_3\}_h-(a_1,\{a_2,a_3\}_h)_h-(-1)^{(|a_3|-1)|a_2|}(\{a_1,a_3\}_h,a_2)_h=\nonumber\\
&&(-1)^{|a_1|+|a_2|-1}(Qn_h'(a_1,a_2,a_3)-n_h'(Qa_1,a_2,a_3)-\nonumber\\
&&(-1)^{|a_1|}n_h'(a_1,Qa_2,a_3)-(-1)^{|a_1|+|a_2|}n_h'(a_1,a_2,Qa_3),\nonumber\\
&&\{\{a_1,a_2\}_h,a_3\}_h-\{a_1,\{a_2,a_3\}_h\}_h+\nonumber\\
&&(-1)^{(|a_1|-1)(|a_2|-1)}\{a_2,\{a_1,a_3\}_h\}_h=0,\nonumber
\end{eqnarray}
where $m_h,m'_h$ are some bilinear operations of degrees $-1$, $-2$ correspondingly and $n_h, n'_h$ are trilinear operations of degree -1, -2 correspondingly.}\\

There exist higher homotopies which turn this homotopy Gerstenhaber algebra into $G_{\infty}$ algebra. This follows from the results of \cite{huangzhao}, \cite{voronov}, \cite{kvz} where it was shown that the symmetrized versions of $(\cdot, \cdot)_h$, $\{,\}_h$, considered on the whole vertex algebra can be continued to the $G_{\infty}$ algebra \cite{tamarkin}. In our case we just need specialization to the conformal weight zero.

One of the central observations of \cite{zeit} was that this  $G_{\infty}$ algebra {\it has a quasiclassical} limit, which can be constructed as follows.
Let $\mathcal{V}\vert_{h=0}=\mathcal{V}^0$ 
(in our example $\mathcal{V}^0=\mathcal{O}(\mathcal{E})$), then consider the subcomplex of  $(\mathcal{F}_h^{\cdot}, Q)$, i.e.  $(\mathcal{F}^{\cdot}, Q)\cong (\mathcal{F}_1^{\cdot}, Q)$, which is:
\begin{eqnarray}\label{subc}
\xymatrixcolsep{30pt}
\xymatrixrowsep{3pt}
\xymatrix{
& \mathcal{V}^0\ar[ddddr] & h\mathcal{V}^0\ar[ddddr]^{\frac{1}{2}h\cdot{\rm div}}& \\
&& \p &&\\
& \bigoplus & \bigoplus & \\
&& {-\frac{1}{2}h\cdot{\rm div}} &&\\
\mathcal{O}_M\ar[uuuur]^{\p} & h\mathcal{O}_M\ar[uuuur]\ar[r]_{i\rm{d}} & h\mathcal{O}_M  & h^2\mathcal{O}_M
}
\end{eqnarray}
The bilinear operations and the operator ${\bf b}$ act on \rf{subc} as follows:
\begin{eqnarray}
&&(\cdot,\cdot)_h: \cF^i\otimes \cF^j\to \cF^{i+j}[h], \quad  \{\cdot,\cdot\}_h: \cF^i\otimes \cF^j\to h\cF_{i+j-1}[h],\\
&&{\bf b}: \cF^i\to h\cF^{i-1}[h], \nonumber
\end{eqnarray}
so that 
\begin{eqnarray}
(\cdot,\cdot)_0= \lim_{h\to 0}(\cdot,\cdot)_h, \quad 
\{\cdot,\cdot\}_0= \lim_{h\to 0}h^{-1}\{\cdot,\cdot\}_h, \quad \mathbf{b}_0=\lim_{h\to 0}h^{-1}\mathbf{b}
\end{eqnarray}
are well defined. The corresponding homotopy Gerstenhaber algebra is 
much less complicated: the corresponding $L_{\infty}$ and $C_{\infty}$ parts are only $L_3$ and $C_3$-algebras. Let us have a look in detail. 
On the level of the vertex algebroid of $\mathcal{O}({\cE})^h$, let us  denote
\begin{eqnarray}\label{vertcour}
&& \lim_{h\to 0}{h}^{-1}[v_1,v_2]=[v_1,v_2]_0, \quad 
\lim_{h\to 0}{h}^{-1}\pi=\pi_0, \\ 
&& \lim_{h\to 0}{h}^{-1}\langle\cdot, \cdot \rangle= \langle\cdot, \cdot \rangle_0.\nonumber 
\end{eqnarray}
Therefore, we can express the bilinear operations $(\cdot, \cdot)_0$ and $\{\cdot,\cdot\}_0$  on the  complex 
\begin{eqnarray}
\xymatrixcolsep{30pt}
\xymatrixrowsep{3pt}
\xymatrix{
& \mathcal{O}(\mathcal{E})\ar[ddddr] & \mathcal{O}(\mathcal{E})\ar[ddddr]^{\frac{1}{2}{\rm div}}& \\
&& d &&\\
& \bigoplus & \bigoplus & \\
&& {-\frac{1}{2}{\rm div}} &&\\
\mathcal{O}_M\ar[uuuur]^{d} & \mathcal{O}_M\ar[uuuur]\ar[r]_{i\rm{d}} & \mathcal{O}_M  & \mathcal{O}_M
}
\end{eqnarray}
 via the following tables:
\begin{center}
$(a_1,a_2)_0$=
\end{center}
\begin{tabular}{|l|c|c|c|c|c|r|}
\hline
\backslashbox{$ a_2$}{$a_1$}&          $u_1$ & $A_1$ & $v_1$ &$\t A_1$ &$\t v_1$& $\t u_1$ \\
\hline
$u_2$ &                  $u_1u_2$    &$A_1u_2$ &$v_1u_2$&
$\t A_1u_2$&$\t v_1 u_2$ &$\t u_1 u_2$\\
&  & $+\pi_0(A_1)u_2$  &  &   && \\
\hline
$A_2$     &  $u_1A_2$ & $-[A_1,A_2]_0-$ &$-v_1A_2$  
&  $\frac{1}{2}\langle\t A_1,A_2\rangle_0$&$-\pi_0(A_2)(\t v_1)$&0\\                    
& & $\frac{1}{2}\langle A_1, A_2\rangle_0$  &   &   && \\
\hline
$v_2$& $u_1\t u_2$ &  $A_1v_2$ & 0 & 0  &$-\t v_1v_2$& 0\\
\hline
$\t A_2$ & $u_1\t A_2$ &$\frac{1}{2}\langle A_1,\t A_2\rangle_0$&  
0 & 0& 0&0\\
\hline
$\t v_2$& $u_1\t u_2$ &   $-\pi_0(A_1)\t v_2$& $-v_1\t v_2$ & 0  &0& 0\\
\hline
$\t u_2$    & $u_1\t u_2$ &   0& 0 & 0  &0& 0\\
\hline
\end{tabular}\\
\vspace{3mm}

\begin{center}
$\{a_1,a_2\}_0$=
\end{center}
\begin{tabular}{|l|c|c|c|c|c|r|}
\hline
\backslashbox{$ a_2$}{$a_1$}&          $u_1$ & $A_1$ & $v_1$ &$\t A_1$ &$\t v_1$& $\t u_1$ \\
\hline
$u_2$ &    0    &$-\pi_0(A_1)u_2$ & 0 & $\pi_0(\t A_1)u_2$& 0 & 0\\
\hline
$A_2$     &  0 & $-[A_1,A_2]_0$ &0   
&  $-[\t A_1,A_2]_0$&$-\pi_0(A_2)\t v_1$&$\pi_0(A_2)\t u_1$\\                    
& &   &  & $-\frac{1}{2}\langle \t A_1, A_2\rangle_0$  && \\
\hline
$v_2$& 0 &  $-\pi_0(A_1)v_2$ & 0 & 0  &$0$& 0\\
\hline
$\t A_2$ & 0 &$-[A_1,\t A_2]_0$&  
0 & $\langle\t A_1,\t A_2\rangle_0$& $-\pi_0(\t A_2)\t v_1$&0\\
\hline
$\t v_2$& 0 &   $-\pi_0(A_1)\t v_2$& 0 & $-\pi_0(\t A_1)\t v_2 $ &0& 0\\
\hline
$\t u_2$    & $-\pi_0(A_1)\t u_2 $&   0& 0 & 0  &0& 0\\
\hline
\end{tabular},\\
\vspace{3mm}

\noindent where $u_i\in\mathcal{F}_h^0$, $(v_i,A_i)\in\mathcal{F}_h^1$, $(\t v_i,\t A_i)\in\mathcal{F}_h^2$, $\t u_i\in \mathcal{F}_h^3$.\\ 

\noindent Let us summarize the results about the quasiclassical limit via Proposition.\\

\noindent{\bf Proposition 3.2.}\cite{zeit} {\it The operations $(\cdot, \cdot)_0$, $\{\cdot,\cdot\}_0$ satisfy the relations \rf{lzrel} so that their symmetrized versions satisfy the relations of $G_{\infty}$ algebra which is the quasiclassical limit of $G_{\infty}$ algebra considered in Proposition 3.1. 
The resulting $C_{\infty}$ and $L_{\infty}$ algebras are reduced to $C_3$ and $L_3$ algebras.  
}\\

The classical limits for the corresponding homotopies $m_h=m_0+O(h)$ and $n_h=n_0+O(h)$ are as follows. The commutativity homotopy $m_0$ is nonzero iff its both arguments belong to $\cF_1$:
\begin{eqnarray}
m_0=-\langle A_1, A_2\rangle_0. 
\end{eqnarray}
The associativity homotopy $n_0$ is nonzero only when all three elements belong to $\cF_1$ or one of the first 
two belongs to $\cF_2$ and the other belong to $\cF_1$:
\begin{eqnarray}
&& n_0(A_1,A_2,A_3)= A_2\langle A_1,A_3\rangle_0-A_1\langle A_2,A_3\rangle_0,\nonumber\\
&& n_0(A_1,\t v,A_2)=n_0(\t v,A_1, A_2)=-\t v \langle A_1, A_2\rangle_0.
\end{eqnarray}

Notice, that in the quasiclassical limit we get rid of all noncovariant terms in the expression for the product and the bracket. This is very close to the classical limit procedure for vertex algebroid. Namely, using \rf{vertcour}, one can obtain vertex algebroid from Courant algebroid. 

The definition of Courant algebroid is as follows (see e.g. \cite{courant}, \cite{bressler}).
A  Courant $\mathcal{O}_M$-algebroid is an $\mathcal{O}_M$-module $\cQ$
equipped with the structure of a Leibniz $\mathbb{C}$-algebra
$[ \cdot, \cdot]_0 : \cQ\otimes_\mathbb{C}\cQ \to \cQ $, 
an $\mathcal{O}_M$-linear map of Leibniz algebras (the anchor map)
$
\pi_0 : \cQ \to \Gamma(TM)
$,
a symmetric $\mathcal{O}_M$-bilinear pairing
$\langle\cdot, \cdot\rangle: \cQ\otimes_{\mathcal{O}_M}\cQ \to \mathcal{O}_M $,  
a derivation
$
\p : \mathcal{O}_M \to \cQ
$,
which satisfy
\begin{eqnarray}
&&\pi\circ\partial =  0, \quad[q_1,fq_2]_0 = f[q_1,q_2] + \pi_0(q_1)(f)q_2,\\
&&\langle [q,q_1],q_2\rangle + \langle q_1,[q,q_2]\rangle  =  \pi_0(q)(\langle q_1, q_2\rangle_0), \quad [q,\partial(f)]_0  =  \partial(\pi_0(q)(f)), \nonumber\\
&&\langle q,\partial(f)\rangle  =  \pi_0(q)(f) \quad [q_1,q_2]_0 + [q_2,q_1]_0  =  \partial(\langle q_1, q_2\rangle_0), \nonumber
\end{eqnarray}
where $f\in\mathcal{O}_M$ and $q,q_1,q_2\in\cQ$. In our case $\cQ\cong\mathcal{O}(\cE)$, $\pi_0$ is just a projection on $\mathcal{O}(TM)$
\begin{eqnarray}
[q_1,q_2]_0=-[q_1,q_2]_D, \quad \langle q_1, q_2\rangle_0=-\langle q_1, q_2\rangle^s, \quad \p=d.
\end{eqnarray}

As we indicated earlier, both $C_{\infty}$ and $L_{\infty}$  parts of $G_{\infty}$ algebra appear to be short. We expect this to happen with all the homotopies, i.e. it is natural to suggest the following.\\

\noindent{\bf Conjecture 3.1.}{\it The $G_{\infty}$ algebra of Proposition 3.2. has only bilinear and trilinear operations, i.e. it is a $G_3$ algebra.}\\

In the following, since we are interested only in the quasiclassical algebra on the complex $(\mathcal{F}^{\cdot},Q)$, we will neglect the 0 subscript for all multilinear operations of this algebra. 

\section{Homotopy Gerstenhaber algebra and Einstein equations}

\noindent{\bf 4.1. BV-subalgebra and a nontrivial example of Einstein equations.}
The homotopy Gerstenhaber algebra we studied in the previous section,  has a subalgebra based on the following complex $(\mathcal{F}^{\cdot}_{sm}, \mathcal{Q})$. 
\begin{eqnarray}
\xymatrixcolsep{30pt}
\xymatrixrowsep{3pt}
\xymatrix{
& \mathcal{O}(T^{(1,0)}M)\ar[ddddr] & \mathcal{O}(T^{(1,0)}M)\ar[ddddr]^{\frac{1}{2}{\rm div}}& \\
&& 0 &&\\
& \bigoplus & \bigoplus & \\
&& {-\frac{1}{2}{\rm div}} &&\\
\mathbb{C}\ar[uuuur]^{0} & \mathbb{C}\ar[uuuur]\ar[r]_{i} & \mathcal{O}_M  & \mathcal{O}_M
}
\end{eqnarray}
It is just a Gerstenhaber algebra (with no higher homotopies), moreover it is a BV algebra, since ${\mathbf b}$ operator also preserves $(\mathcal{F}^{\cdot}_{sm}, Q)$. Therefore, we have the following Proposition.\\

\noindent{\bf Proposition 4.1.} {\it Bilinear operations $(\cdot,\cdot)$, $\{\cdot,\cdot\}$ together  with operator $\mathbf{b}$ generate the structure of BV algebra on $(\mathcal{F}^{\cdot}_{sm}, Q)$.}\\
 
Let us consider the ${\infty}$-jet version of the complex $(\mathcal{F}^{\cdot}_{sm}, Q)$: we substitute $\mathcal{O}_M$, $\mathcal{O}(T^{(0,1)}(M))$ by $J^{\infty}(\mathcal{O}_M)$, $J^{\infty}(\mathcal{O}(T^{(0,1)}(M))$. We denote the  resulting complex as $(\mathcal{F}^{\cdot}_{sm, \infty}, Q)$. 
Then the completed tensor product 
\begin{eqnarray}
{\bf F}^{\cdot}_{sm, \infty}=\mathcal{F}^{\cdot}_{sm,\infty}\hat{\otimes}\bar {\mathcal{F}}^{\cdot}_{sm, \infty} 
\end{eqnarray}
where $(\bar{\mathcal{F}}^{\cdot}_{sm, \infty}, \bar{Q})$ is the  complex obtained from $(\mathcal{F}^{\cdot}_{sm, \infty}, Q)$ by complex conjugation. Complex 
$({\bf F}^{\cdot}_{sm, \infty}, \mathcal{Q})$, where $\mathcal{Q}=Q+\bar{Q}$, 
is the jet version of the complex $({\bf F}_{sm}^{\cdot}, \mathcal{Q})$, such 
that e.g. ${\bf F}_{sm}^2=\Gamma(T^{(1,0)}M\otimes T^{(0,1)}M)\oplus (\bar{\mathcal{O}}(T^{(0,1)}M)\oplus \mathcal{O}(T^{(1,0)}M))^{\oplus 2}\oplus
\mathcal{O}_M\oplus \bar{\mathcal{O}}_M\oplus \mathbb{C}$. 
Clearly, the complex $({\bf F}_{sm}^{\cdot}, \mathcal{Q})$ carries a 
structure of BV algebra inherited from $(\mathcal{F}^{\cdot}_{sm, \infty}, {Q})$ and its complex conjugation, so that 
\begin{eqnarray}\label{brack2}
&&(-1)^{|a_1|}\{a_1,a_2\}=\\
&&\mathbf{b^-}(a_1,a_2)-(\mathbf{b^-}a_1,a_2)-(-1)^{|a_1|}(a_1\mathbf{b^-}a_2)\nonumber,
\end{eqnarray}
where $\mathbf{b^-}=\mathbf{b}-\bar{\mathbf{b}}$. Note, that the elements closed under $\mathbf{b}^-$ form a subalgebra in the differential graded algebra (DGLA), generated by $Q, \{\cdot, \cdot\}$. 
It turns out that the Maurer-Cartan equations of this DGLA and their symmetries have a very interesting meaning. To describe them, let us define some extra algebraic operations for convenience.

Let $g,h \in \Gamma(T^{(1,0)}M\otimes T^{(0,1)}M)$ so that their components are  
$g^{i\b j}\p_i\otimes\p_{\b j}, h^{i\b j}\p_i\otimes\p_{\b j}$.  
Then one can define symmetric bilinear operation \cite{lmz}, \cite{zeit2}:
\begin{eqnarray}
&&[[,]]:\Gamma(T^{(1,0)}M\otimes T^{(0,1)}M)\otimes \Gamma(T^{(1,0)}M\otimes T^{(0,1)}M)\to\\  
&& \Gamma(T^{(1,0)}M\otimes T^{(0,1)}M)\nonumber
\end{eqnarray}
written in components as follows:
\begin{eqnarray}
&&[[g,h]]^{k\b l}\equiv\\
&&(g^{i\b j}\p_i\p_{\b j}h^{k\b l}+h^{i\b j}\p_i\p_{\b j}g^{k\b l}-\p_ig^{k\b j}\p_{\b j}h^{i\b l}-
\p_ih^{k\b j}\p_{\b j}g^{i\b l})\nonumber
\end{eqnarray}
and looks much less complicated in the jet notation (see section 2). Namely, if $\tilde{\xi}\tilde{\eta}\in J^{\infty}(\mathcal{O}(T^{(1,0)}M)\otimes J^{\infty}(\mathcal{O}(T^{(0,1)}M)$, so that $\xi=\sum_I v^I\otimes\b{v}^I$, $\eta=\sum_J w^J\otimes\b{w}^J$, where $v^I, w^J\in J^{\infty}(\mathcal{O}(T^{(1,0)}M)$, $\b{v}^I, \b{w}^J\in J^{\infty}(\mathcal{O}(T^{(0,1)}M)$, then 
\begin{eqnarray}
[[\xi,\eta]]=\sum_{I,J}[v^I,w^J]\otimes[\b v^I,\b w^J].
\end{eqnarray} 

As noted in \cite{lmz},\cite{zeit2}, if the bilinear tensor 
$g$ is such that one can associate a  
K\"ahler metrics to it, then the Ricci tensor $R^{i\b j}$ associated with such metric tensor is proportional to 
$[[g,g]]$, more precisely 
\begin{eqnarray}
R^{i\b j}(g)=\frac{1}{2}[[g,g]]^{i\b j}.
\end{eqnarray}

If the complex manifold $M$ has a volume form 
$\Omega$, such that in local coordinates 
$\Omega=e^f dX^1\dots\wedge dX^n\wedge dX^{\bar 1}\wedge\dots dX^{\bar n}$. 
Let us denote the volume form which determines the differential $\mathcal{Q}$ as $\Omega'$, so that $f=-2\Phi_0'$, then $\Phi_0'$ has to be locally a sum of holomorphic and antiholomorphic functions, i.e. it satisfies equation  
$\p_i{\p}_{\b j}\Phi_0'=0$. 

We will refer to the vector field $div_{\Omega}g$ such that 
$(div_{\Omega}g)^{\b j}=\p_ig^{i\b j}+\p_ifg^{i\b j}$, $(div_{\Omega}g)^{i}=\p_{\b j}g^{i\b j}+\p_{\b j}fg^{i\b j}$ as the divergence of bivector field $g$ with respect to the volume form $\Omega$.

Now let the Maurer-Cartan element, closed under $\mathbf{b}^-$, namely the element of $
\Gamma(T^{(1,0)}M\otimes T^{(0,1)}M)\oplus \mathcal{O}(T^{(0,1)}M)\oplus \mathcal{O}(T^{(1,0)}M)\oplus
\mathcal{O}_M\oplus \bar{\mathcal{O}}_M$ be defined by its components in the  direct sum, i.e. as $(g, \b v, v, \phi, \b \phi)$.

Then the following Theorem holds, which can be proven by direct calculation.\\

\noindent{\bf Theorem 4.1a.} {\it 
The Maurer-Cartan equation for the differential graded Lie algebra on ${\bf F}^{\cdot}_{sm}\vert_{\mathbf{b^-}=0}$ generated by $\mathcal{Q}$ and $\{\cdot, \cdot\}$ imposes the following system of equations on $g, \phi, \b \phi$ ($\b v, v$ turn out to be auxilliary variables): }

1). {\it Vector field $div_{\Omega}g$, where $\Omega=\Omega'e^{-2\phi+2\b \phi}$ is determined by $f\equiv-2\Phi_0=-2(\Phi_0'+\phi-\b \phi)$ and $\p_i{\p}_{\b j}\Phi_0=0$, is such that its $\Gamma(T^{(1,0)} M)$, $\Gamma(T^{(0,1)} M)$ components are correspondingly holomorphic and antiholomorphic.}

2). {\it Bivector field $g\in \Gamma(T'M\otimes T''M)$ obeys the following equation: 
\begin{eqnarray}\label{bil}
[[g,g]]+\mathcal{L}_{div_{\Omega}(g)}g=0, 
\end{eqnarray} 
where $\mathcal{L}_{div_{\Omega}(g)}$ is a Lie derivative with respect to 
the corresponding vector fields.}\\

3). {\it $div_{\Omega}div_{\Omega}(g)=0$.\\

The infinitesimal symmetries of the Maurer-Cartan equation coincide with the holomorphic coordinate transformations of the volume form and tensor $\{g^{i\b j}\}$. }\\

The constraints 1), 2), 3) coincide with the equations studied in \cite{zeit2},  where it was shown that they are equivalent to Einstein equations, i.e. the following statement is valid.\\ 

\noindent{\bf Theorem 4.1b.} 
{\it If tensor $\{g^{i\b j}\}$ parametrises Hermitean metric, then 
the conditions 1), 2), 3) on $g$ and $\Phi_0$ from Theorem 4.1a are equivalent to Einstein equations 
\begin{eqnarray}\label{components}
&&R^{\mu\nu}={1\over 4} H^{\mu\lambda\rho}H^{\nu}_{\lambda\rho}-2\nabla^{\mu}
\nabla^{\nu}\Phi,\\
&&\nabla_{\mu}H^{\mu\nu\rho}-2(\nabla_{\lambda}\Phi)H^{\lambda\nu\rho}=0,
\nonumber\\
&&4(\nabla_{\mu}\Phi)^2-4\nabla_{\mu}\nabla^{\mu}\Phi+
R+{1\over 12} H_{\mu\nu\rho}H^{\mu\nu\rho}=0,\nonumber
\end{eqnarray}
where $H=dB$ is a 3-form, so that metric $G$, 2-form $B$ and the dilaton field $\Phi\in \mathcal{C}(M)$ are expressed as follows: 
\begin{eqnarray}\label{constr}
&&G_{i\bar{k}}=g_{i\bar{k}}, \quad B_{i\bar{k}}=-g_{i\bar{k}}, \quad \Phi=\log\sqrt{g}+\Phi_0,\\
&&G_{ik}=G_{\b i \b k}=G_{ik}=G_{\b i \b k}=0,\nonumber
\end{eqnarray}
where by $g$ under the square root we denote the determinant of $\{g_{i\b j}\}$. 
In other words, \rf{components} are equivalent to the following system:
\begin{eqnarray}
&&\p_i\p_{\bar{k}}\Phi_0=0,\quad \p_{\bar{p}}d^{\Phi_0}_{\bar{l}}g^{\bar{l}k}=0, 
\quad \p_{p}d^{\Phi_0}_{l}g^{\bar{k}l}=0,\nonumber\\
&&2g^{r\bar{l}}\p_r\p_{\bar{l}}g^{i\bar{k}}-2\p_r g^{i\bar{p}}\p_{\bar{p}}g^{r\bar{k}}-
g^{i\bar{l}}\p_{\bar{l}}d^{\Phi_0}_sg^{s\bar{k}}-g^{r\bar{k}}\p_r d^{\Phi_0}_{\bar{j}}g^{\bar{j}i}+\nonumber\\
&&\p_rg^{i\bar{k}}d^{\Phi_0}_{\bar{j}}g^{\bar{j}r}+\p_{\bar{p}}g^{\bar{k}i}d^{\Phi_0}_n g^{n\bar{p}}=0,\nonumber\\
&&d^{\Phi_0}_{i}d^{\Phi_0}_{\bar{j}}g^{\bar{j}i}=0,
\end{eqnarray}
where $d^{\Phi_0}_{i}g^{\bar{j}i}\equiv (\p_i -2\p_i\Phi_0) g^{i\bar{j}}$, $d^{\Phi_0}_{\bar{j}}g^{\bar{j}i}\equiv (\p_{\b j}-2\p_{\b j}\Phi_0) g^{i\bar{j}}$, which is the component form 
of the conditions 1), 2) and 3). }\\

\noindent{\bf 4.2. Physical motivation for the main conjectures.} 
In this subsection, we give the physics motivation for the generalization of the  result of Section 4.1: namely, we want to extend Theorem 4.1. to the case of full complex $\mathcal{F}^{\cdot}$. 
For more details we refer the reader to the paper \cite{lmz}, \cite{zeit2}. In Section 2, we considered the equivalence of two action functionals $S_{fo}$ and $S_{so}$. 
On the quantum level the object of primary interest is the path integral 
\begin{eqnarray}
\int [dp][d\b p][dX][d\b X]e^{-S}.
\end{eqnarray}
In the case of $S=S_0$ the quantum theory corresponding to this path integral, is described by the gerbes of chiral differential operators (and locally just by vertex algebras), as it was already mentioned in Section 3. However, this action should be modified to accomodate the holomorphic volume form on $M$, otherwise the Virasoro element in the corresponding vertex algebras wouldn't be globally defined.  
On the level of action functionals, one has to add an extra term to $S_0$, namely  
\begin{eqnarray}
S_0\to S_0+\int_{\Sigma}\sqrt{\gamma}R^{(2)}(\gamma)\phi(X), 
\end{eqnarray}
where $e^{-2\phi}$ is the density for the volume form on $M$, so that $\p_i\p_{\b j}\phi=0$, $\gamma$ is a metric on $\Sigma$ and $R^{(2)}(\gamma)$ is its curvature. Let us add a similar term to its perturbed version $S_{fo}$, i.e. $\int_{\Sigma}\sqrt{\gamma}R^{(2)}(\gamma)\Phi_0(X)$ with no restrictions on $\Phi_0$, and the resulting action will be denoted as $S^{\Phi_0}_{fo}$. 
We will call $\Phi_0$ a $normalized$ $dilaton$ $field$.

The integration over $p,\b p$ leads to the following (see \cite{lmz}):
\begin{eqnarray}
\int [dp][d\b p][dX][d\b X]e^{-S^{\Phi_0}_{fo}}=
\int [dX][d\b X]e^{-S_{so}+\int R^{(2)}(\gamma)(\Phi_0(X)+\sqrt{g})}.
\end{eqnarray}
This heuristic formula gives the proper correspondence between the first- and second-order actions on the quantum level. 
The relation of those to Einstein equations is as follows. 
Analysing the path integral in the right hand side involves regularization procedure which leads to the broken scale (and conformal) invariance. 
In order to make the model conformally invariant, one has to impose a sequence of constraints, appearing for the vanishing of the $\beta$-function (see \cite{eeq}, \cite{fts}, \cite{polbook}). The $\beta$-function is the function depending on $G$-metric, 2-form $B$, dilaton 
\begin{eqnarray}\label{dil}
\Phi=\Phi_0+\sqrt{g}
\end{eqnarray} 
and parameter $h$. At the zeroth order in $h$ vanishing of the $\beta$-function leads to Einstein equations \rf{components}. Vanishing of the coefficients of higher powers in $h$ lead to the equations involving higher number of derivatives and higher powers in Ricci curvature. It was noted (see e.g. \cite{pol}) that the linearized form of Einstein equations and their symmetries can be obtained as the closedness condition for the elements of degree 2 in  the semi-infinite (BRST) complex associated to the Virasoro module corresponding to the conformal field theory described by $S_{so}$ with the flat metric. One of the statements of String Field Theory is that the full conformal invariance conditions can be obtained from Maurer-Cartan equations for some $L_{\infty}$-algebra on BRST complex \cite{zwiebach}, so that the full metric, B-field and dilaton can be restored from the Maurer-Cartan element. The symmetries of the Maurer-Cartan equations correspond to the $h$-corrected diffeomorphism symmetries and the exact shifts of the antisymmetric tensor $B$.

The complex corresponding to the flat 
metric does not have any simple algebraic structure on it (because it is not a vertex algebra) and it is complicated to construct such algebraic operations explicitly.  
On the contrary, for the first-order model we start from the vertex 
algebra and we have related $G_{\infty}$ structure due to \cite{lz}, \cite{kvz}, \cite{huangzhao}, \cite{voronov}. Using the results of \cite{zeit}, we are able to reduce it to much smaller complex and find the quasiclassical limit.  
We claim 
that the Maurer-Cartan equation corresponding to its $L_{\infty}$-subalgebra of the  quasiclassical limit of this $G_{\infty}$ algebra reproduces Einstein equations and their symmetries, where the metric, 2-form $B$ and the dilaton $\Phi$ are expressed by means of \rf{GB}, 
\rf{dil}. 
In subsection 4.1., we obtained this correspondence in the case when only one of the perturbing terms was present in $S_{fo}$, namely $\langle g,p\wedge\b p\rangle $.
In the next subsection, we extend the statement of Theorem 4.1 to the case of general $S_{fo}$.\\

\noindent{\bf 4.3. Main Conjectures.} Following the ideas of Section 4.1, we want to repeat the construction in the case of the complex $(\mathcal{F}^{\cdot}, \mathcal{Q})$. Namely, we consider its jet version $(\mathcal{F}^{\cdot}_{\infty}, \mathcal{Q})$
and its complex conjugate $(\bar{\mathcal{F}}^{\cdot}_{\infty}, \mathcal{Q})$, so that
\begin{eqnarray}
{\bf F}^{\cdot}_{ \infty}=\mathcal{F}^{\cdot}_{\infty}\hat{\otimes}\bar {\mathcal{F}}^{\cdot}_{\infty}. 
\end{eqnarray}
It is the jet version of the complex $({\bf F}^{\cdot}, \mathcal{Q})$, such that e.g. the subspace of degree 1 is as follows:  ${\bf F}^{1}\cong\Gamma(E)\oplus\mathcal{C}(M)\oplus\mathcal{C}(M)$. As in the section 4.1, the  divergence operator which determines $\mathcal{Q}$-operator, is based on the volume form, given in the local coordinates as $e^{-2\Phi_0'(X)}dX^1\dots dX^n\wedge dX^{\b 1}\dots dX^{\b n}$, so that $\p_i\p_{\b j}\Phi_0'=0$.

We can give ${\bf F}^{\cdot}$ the structure of the homotopy Leibniz bracket by means of formula  which is the same as in Section 4.1:  
\begin{eqnarray}
&&(-1)^{|a_1|}\{a_1,a_2\}=\nonumber\\
&&\mathbf{b^-}(a_1,a_2)-(\mathbf{b^-}a_1,a_2)-(-1)^{|a_1|}(a_1\mathbf{b^-}a_2),
\end{eqnarray}
however now we have higher homotopies. 
We also note that ${\bf F}^{\cdot}\vert_{\mathbf{b}^-=0}\equiv {\bf F}^{\cdot}_-$ is invariant under $\{\cdot,\cdot\}$. 
Let us formulate the first part of the main conjecture.\\

\noindent{\bf Conjecture 4.1a.} {\it The structure of homotopy Gerstenhaber algebra on 
${\bf F}^{\cdot}$ can be extended to $G_{\infty}$-algebra, so that the subcomplex 
${\bf F}^{\cdot}_-$ is invariant under $L_{\infty}$ operations. }\\

Let us focus on the subcomplex $(\mathbf{F}^{\cdot}_-, \mathcal{Q})$. The space of Maurer-Cartan elements, i.e. the subspace of the elements of degree 2 is:
\begin{eqnarray}
\mathbf{F}_-^2\cong \Gamma(\mathcal{E}\otimes\bar{\mathcal{E}})\oplus \Gamma(E)\oplus \mathcal{C}(M)\oplus\mathcal{C}(M).
\end{eqnarray}
The elements of this space are defined by means of the components from the direct sum above, i.e. $\Psi=(\mathbb{M}, \eta, \phi,\bar{\phi})$. We will denote the difference $\phi-\bar{\phi}\equiv\Phi_0''$ and $\Phi_0\equiv\Phi_0'+\Phi_0''$. 
Let us formulate the second part of the main conjecture.
\\

\noindent{\bf Conjecture 4.1b.} {\it Let $\Psi=(\mathbb{M}, \eta, \phi,\bar{\phi})$ be the solution of the generalized Maurer-Cartan (GMC) equation for $L_{\infty}$-algebra on ${\bf F}^{\cdot}_-$, so that 
\begin{equation}
\mathbb{M}=\begin{pmatrix} g & \mu \\ 
\bar{\mu} & b \end{pmatrix},
\end{equation} 
Then the $\eta$-component is auxilliary and is expressed via $\mathbb{M}$ and $\phi, \bar{\phi}$. If $\{g^{i\b j}\}$ is invertible, then $G, B$ obtained from $\mathbb{M}$ via \rf{GB} together with $\Phi=\Phi_0+\sqrt{g}$, where $g$ is the determinant of $\{g_{i\b j}\}$, satisfy the Einstein equations \rf{components}.}\\
 
The space of infinitesimal symmetry generators of GMC equation, i.e. $\mathbf{F}^1$ is given by  
\begin{eqnarray}
\mathbf{F}_-^1\cong\Gamma(E)\oplus \mathcal{C}(M),
\end{eqnarray}
so that any element can be written in components as $\Lambda=(\xi, f)$.

The third part of the conjecture concerns the question how $\xi, f$ are related to 
$\alpha\in \Gamma(E)$ in the transformation formula 
\begin{eqnarray}\label{mtransf}
\mathbb{M}\to \mathbb{M}-D\alpha+ \phi_1(\alpha,\mathbb{M})+\phi_2(\alpha, \mathbb{M},\mathbb{M})
\end{eqnarray} 
from Section 1. 
First, to justify the statement of Conjecture 4.1b., we prove the following Proposition. \\

\noindent{\bf Proposition 4.2.} {\it 
Let $\Lambda=(\xi, f)\in \mathbf{F}_-^1$ be the generator of the  infinitesimal transformation of the solution of GMC equation. Then after the substitution   
$\xi=\alpha+\frac{1}{2}\mathbb{M}\cdot\alpha$ (where $\mathbb{M}$ is considered as an element of $End(\Gamma(E))$) the transformation of $\mathbb{M}$-component of the solution coincides with \rf{mtransf} up to the second order in $\mathbb{M}$.}\\

\noindent {\bf Proof.} At first, we show that the expression 
\begin{eqnarray}\label{mtransfch}
\Psi\to \Psi+\mathcal{Q}\Lambda'-\{\Lambda', \Psi\}+\frac{1}{2}\{\Lambda',\Psi, \Psi\} 
\end{eqnarray}
gives the formula \rf{mtransf}, where $\{\cdot, \cdot, \cdot\}$ is the homotopy for the Leibniz identity for $[\cdot, \cdot]$, $\Lambda'=(\alpha,s)$, $\Psi=(\mathbb{M}, \eta, \phi,\bar{\phi})$. 
It is easy to check the corresspondence between \rf{mtransf} and \rf{mtransfch}  for the 0th and 1st order in $\mathbb{M}$.  To prove Proposition 4.2, we just need to check the term corresponding to trilinear operation. Let us return to the jet level, i.e. we assume that
\begin{eqnarray}
&&\alpha=\sum_Jf^J\otimes {\bar{b}}^J+\sum_Kb^K\otimes {\bar{f}}^K,\nonumber\\ 
&&\mathbb{M}=\sum_I a^I\otimes \bar{a}^I.
\end{eqnarray}
Then the only terms contributing to the relevant part of $\{\Lambda',\Psi, \Psi\}$ are as follows:
\begin{eqnarray}
&&-\sum_{I,J,K}(m_0'(b^I, a^K), a^J)\otimes (\{\bar{a}^J,\bar{f}^I\}, {\bar{a}}^K)-\nonumber\\
&&\sum_{I,J,K}
 (\{a^J,f^I\}, {a}^K)\otimes (m_0' ({\bar{b}}^I, {\bar{a}}^K), {\bar{a}}^J).
\end{eqnarray}
We see that modulo the necessary coefficient this coincides with the trilinear operation $\phi_2$ \rf{trialgop}. 
The statement of the Proposition 4.2 can be obtained from the antisymmetrization of 
$\{\cdot,\cdot \}$, and therefore of $\{\cdot,\cdot, \cdot \}$, 
so that the formula
\begin{eqnarray}\label{sym}
\Psi\to \Psi+\mathcal{Q}\Lambda+\{\Psi,\Lambda\}^{asymm}+\frac{1}{2}\{\Psi, \Psi, \Lambda\}^{asymm}+\dots 
\end{eqnarray}
corresponds to \rf{mtransf} if $\Lambda=(\alpha+\frac{1}{2}\mathbb{M}\cdot \alpha, s).$ 
\hfill$\blacksquare$

We note, that the symmetry generated by $f$-part of $\mathbf{F}^1$ element does not affect metric B-field or dilaton. It is easy to check that on  the level of 0th order in $\mathbb{M}$: the symmetry transformation corresponds to the shift of $\phi $ and $\b \phi$  by $f$. One can check, similar to Proposition 4.2, that this symmetry remains redundant for the first and second order. We claim that these statements are exact, namely the following Conjecture is true.\\

\noindent{\bf Conjecture 4.1c.} {\it Let $\Lambda=(\xi,f)\in \mathbb{F}^1_-$, be the generator of the infinitesimal symmetries of  GMC equation \rf{sym}. The corresponding transformation of $\mathbb{M}$-component of the solution of GMC coincide with \rf{mtransf} if  $\xi=\alpha+\frac{1}{2}\mathbb{M}\cdot\alpha$. 
Under conditions of Conjecture 4.1b these transformations  reproduce infinitesimal diffeomorphism transformations and shifts of B-field by exact 2-form, which are the symmetries of equations \rf{components}.}

\section{Appendix}
Here we give explicitly the formulas for the transformations of $\mathbb{M}$ (see Section 2): 
\begin{eqnarray}
\mathbb{M}\to \mathbb{M}-D\alpha+ \phi_1(\alpha,\mathbb{M})+\phi_2(\alpha, \mathbb{M},\mathbb{M})
\end{eqnarray}
of the matrix elements of Beltram-Courant differential 
\begin{equation}
\mathbb{M}=\begin{pmatrix} g & \mu \\ 
\bar{\mu} & b \end{pmatrix},
\end{equation}
where $g\in \Gamma(T^{(1,0)}M\otimes T^{(0,1)}M)$, $\mu\in \Gamma(T^{(1,0)}M\otimes {T^*}^{(0,1)}M)$, $\b \mu\in \Gamma({T^*}^{(1,0)}M\otimes T^{(0,1)}M)$, $b\in \Gamma({T^*}^{(1,0)}M\otimes {T^*}^{(0,1)}M)$. The explicit form of the transformations in components is (with the notations from Section 2):
\begin{eqnarray}
&&g^{i{\bar j}} \rightarrow g^{i{\bar j}} + v^k\p_k g^{i{\bar j}} + v^{\bar{k}}\p_{\bar{k}} g^{i{\bar j}}
- g^{i{\bar k}}\d_{\bar k}v^{\bar j}
- g^{k{\bar j}}\d_kv^i+\\
&&g^{i{\bar k}}{\bar\mu}^{\bar j}_k\d_{\bar k}v^k+g^{k{\bar j}}\mu^{i}_{\bar{k}}\d_{k}v^{\bar{k}},\nonumber\\
 &&\mu^{i}_{\bar{j}} \rightarrow 
 \mu^{i}_{\bar{j}} -
\p_{\bar{j}}v^i + v^{k}\p_k\mu^{i}_{\bar{j}} +
v^{\bar{k}}\p_{\bar{k}}\mu^{i}_{\bar{j}}+
\mu^{i}_{\bar{k}}\p_{\bar{j}}v^{\bar{k}} -
\mu^k_{\bar{j}}\p_kv^i+ \nonumber\\
&&\mu^i_{\bar{l}}\mu^k_{\bar{j}}\p_k v^{\bar{l}}+g^{i\b k}b_{l\b j}\p_{\b k}v^l,\nonumber\\
&&b_{i{\bar j}} \rightarrow b_{i{\bar j}} + v^k\p_k b_{i{\bar j}} + v^{\bar{k}}\p_{\bar{k}} b_{i{\bar j}}
+ b_{i{\bar k}}\p_{\bar{j}}v^{\bar{k}}+b_{l{\bar j}}\p_i v^l+\nonumber\\
&&b_{i{\bar k}}\mu^{k}_{\bar j}\p_kv^{\bar{k}}
+b_{l{\bar j}}{\bar\mu}^{\bar k}_i\p_{\bar k}v^l,\nonumber
\end{eqnarray}
 
\begin{eqnarray}
&&g^{i\b j}\to g^{i\b j}+g^{i\b k}(\p_l\omega_{\b k}-\p_{\b k}\omega_l)g^{l\b j},\\
&&\mu^i_{\b l}\to \mu^i_{\b l}+g^{i\b j}(\p_{\b l}\omega_{\b j}-\p_{\b j}\omega_{\b l})+
\mu^r_{\b l}g^{i\b k}(\p_r\omega_{\b k}-\p_{\b k}\omega_r), \nonumber\\
&&b_{i\bar{j}}\to  b_{i\bar{j}}+\p_{\bar j}\omega_i-\p_i\omega_{\bar j}+\mu^i_{\bar j}(\p_i\omega_k-\p_k\omega_i)+\bar{\mu}^{\bar s} _i(\p_{\bar j}\omega_{\bar s}-
\p_{\bar s}\omega_{\bar j})+\nonumber\\
&&{\bar \mu}^{\bar i}_j\mu_{\bar k}^s(\p_s\omega_{\bar i}-\p_{\bar i}\omega_s).\nonumber
\end{eqnarray}

\end{document}